\documentclass[12pt]{amsart}
\usepackage{amsfonts,latexsym}

\newtheorem{theorem}{Theorem}
\newtheorem{lemma}[theorem]{Lemma}
\newtheorem{thm}[theorem]{Theorem}

\theoremstyle{definition}

\theoremstyle{remark}
\newtheorem{remark}[theorem]{Remark}

\title[Combinatorial principles, compactness of spaces]
{Combinatorial and model-theoretical principles related to 
regularity of ultrafilters and compactness of topological spaces. II.}

\author[]{Paolo Lipparini} 
\address{Dipartimento di Matematica\\
Viale della Ricerca Scientifica\\
II Universit\`a Romanaccia (Tor Vergata)\\
I-00133 ROME 
ITALY
}
\urladdr{http://www.mat.uniroma2.it/\textasciitilde lipparin}
\thanks{The author has received support from MPI and GNSAGA.
We wish to expressed our gratitude to X. Caicedo for stimulating discussions and correspondence} 


\begin{document} 

\begin{abstract} 
We find many conditions equivalent to the model-theoretical property 
$\lambda \stackrel{\kappa}{\Rightarrow} \mu$
introduced in \cite{bumi}.
 Our conditions involve uniformity of ultrafilters, compactness
properties of products of topological spaces and the existence  of
 certain
 infinite matrices. 
\end{abstract}

\maketitle




See Part I \cite{parti} or \cite{CN,CK,KM,KV,EGT} for unexplained 
notation.

According to \cite{bumi}, if
$ \lambda \geq \mu$ are infinite regular cardinals, and $ \kappa $
is a cardinal,  
$\lambda \stackrel{\kappa}{\Rightarrow} \mu$
 means that the model $ \langle \lambda, <, \gamma \rangle _{ \gamma < \lambda }  $ has an
expansion ${\mathfrak A}$ in a language with at most $ \kappa $ new symbols such that whenever
$\mathfrak B \equiv \mathfrak A$ and $ \mathfrak B$ has an element $x$ such that 
$ {\mathfrak B} \models \gamma < x $ for every $ \gamma < \lambda $, 
then   
$ \mathfrak B$ has an element $y$ such that 
$ {\mathfrak B} \models \alpha  < y < \mu $ for every $ \alpha  < \mu $.

An ultrafilter $D$ over $ \lambda $ is said to be 
uniform if and only if every member of $D$ 
has cardinality $ \lambda $. If $ \lambda $ is a regular cardinal, then it is
obvious that an ultrafilter $D$ is uniform over $ \lambda $ if and only if 
the interval $[0, \gamma] \not \in D$, for every $ \gamma < \lambda $, if and only if 
the interval $(\gamma , \lambda )$ is in $D$, for every $ \gamma < \lambda $. 

Thus, if  $D$ is an ultrafilter
 over some regular cardinal $ \lambda $, and if $Id_D$
denotes the $D$-class  of the identity
function  on $ \lambda $, then $D$ is uniform
over $ \lambda $ if and only if   
in the model
  $ {\mathfrak C} = \prod_D \mathfrak A$  we have that
$ d( \gamma) < Id_D$ for every $ \gamma < \lambda $.
Here, $d$  denotes the elementary embedding.

If $D$ is an ultrafilter over $I$, and $f:I\to J$, then 
$f(D)$ is the ultrafilter over $J$ defined by:
$Y \in f(D)$ if and only if $f ^{-1} (Y) \in D $.  

If $ \kappa  , \lambda  $ are infinite cardinals, a topological space is
said to be  {\em $[ \kappa  , \lambda  ]$-compact}
if and only if
every open cover by at most $ \lambda $  sets  has a subcover by less than 
$\kappa $ sets.
No separation axiom is needed to prove the results of the present paper.

\begin{thm}\label{lmkprod2} 
Suppose that  $ \lambda \geq \mu$  are infinite regular cardinals, 
and $\kappa \geq \lambda $ is an infinite cardinal.  
Then the following conditions are equivalent.

\smallskip

(a) $\lambda \stackrel{\kappa}{\Rightarrow} \mu$ holds.

\smallskip

(b) There are $ \kappa $ functions $ (f_ \beta ) _{ \beta < \kappa } $
from $ \lambda $ to $\mu$ such that whenever $D$ is an ultrafilter
uniform over $ \lambda $ then there exists some $ \beta < \kappa $
such that $f_ \beta (D)$ is uniform over $ \mu$.

\smallskip

(b$'$) There are $ \kappa $ functions $ (f_ \beta ) _{ \beta < \kappa } $
from $ \lambda $ to $\mu$ for which the following holds:
for every function $g: \kappa \to \mu$ there exists some finite
set $F \subseteq \kappa $ such that 
$ \left| \bigcap _{\beta \in F} f_\beta  ^{-1}([0, g(\beta ))) \right| < \lambda $.

\smallskip

(c) There is a family $ (B_{ \alpha , \beta }) _{ \alpha<\mu , \beta<\kappa}  $ 
of subsets of $ \lambda $ such that:

(i) For every $ \beta<\kappa$, $\bigcup _{ \alpha<\mu } B_{ \alpha , \beta  } = \lambda$;

(ii) For every $ \beta<\kappa$ and $ \alpha \leq \alpha ' < \mu  $, 
$ B_{ \alpha , \beta } \subseteq B_{ \alpha' , \beta }$;

(iii) For every function $g : \kappa  \to \mu $ there exists a finite subset
$F \subseteq \kappa  $ such that 
$|\bigcap _{\beta \in F} B_{ g( \beta) , \beta }| < \lambda $.   

\smallskip

(d)
Whenever
$(X_ \beta ) _{ \beta < \kappa }$ is a family of topological spaces
such that  no $X_ \beta $ is 
$[ \mu, \mu]$-compact,
then $X=\prod_{ \beta < \kappa } X_ \beta $
 is not  $[ \lambda , \lambda ]$-compact.

\smallskip

(e) The topological space $ \mu^ \kappa $ is not
$[ \lambda , \lambda ]$-compact, where $ \mu$
is endowed with the topology whose open sets are the
intervals $ [0, \alpha) $ ($ \alpha \leq \mu$), and 
$ \mu^ \kappa $ is endowed with the Tychonoff topology. 
\end{thm} 
 
\begin{remark}\label{future}
An analogue of Theorem \ref{lmkprod2}
holds for the more general notion
$(\lambda, \mu) \stackrel{\kappa}{\Rightarrow} (\lambda', \mu')$
introduced in \cite{easter} (see also \cite[Section 0]{arxiv}).
Details shall be presented elsewhere.
For this more general notion, the equivalence of conditions 
analogue to
(a) and (b) above has been stated in 
\cite{abst}. 
There we also stated the analogue of 
(b) $ \Rightarrow $ (d).
\end{remark}

\begin{proof}
(a) $\Rightarrow$ (b).
Let ${\mathfrak A}$
be an expansion of 
$ \langle \lambda, <, \gamma \rangle _{ \gamma < \lambda }  $
witnessing
$\lambda \stackrel{\kappa}{\Rightarrow} \mu$.

Without loss of generality we can assume that 
${\mathfrak A}$ has Skolem functions (see \cite[Section 3.3]{CK}).
Indeed, since $ \kappa \geq \lambda $, adding Skolem functions
to ${\mathfrak A}$ involves adding at most $ \kappa $ new symbols.

Consider the set of all functions $f: \lambda \to \mu$ 
which are definable in ${\mathfrak A}$. Enumerate them as 
$ (f_ \beta ) _{ \beta < \kappa } $. We are going to show that 
these functions witness (b).

Indeed, let $D$ be an ultrafilter
uniform over $ \lambda $. 
Consider the $D$-class $Id_D$ of the identity
function  on $ \lambda $. 
Since $D$ is
uniform over $ \lambda $, in the model
  $ {\mathfrak C} = \prod_D \mathfrak A$  we have that
$ d( \gamma) < Id_D$ for every $ \gamma < \lambda $, 
where $d$  denotes the elementary embedding.
Let $ {\mathfrak B} $ be the Skolem hull
of $Id_D$ in ${\mathfrak C}$. 
By 
\L o\v s Theorem,
$ {\mathfrak C} \equiv {\mathfrak A} $.
Since 
$\mathfrak A $ has Skolem functions, 
$ {\mathfrak B} \equiv {\mathfrak C}$ \cite[Proposition 3.3.2]{CK}. By transitivity,
$ {\mathfrak B} \equiv {\mathfrak A} $.

Since ${\mathfrak A}$
witnesses
$\lambda \stackrel{\kappa}{\Rightarrow} \mu$,
then   
$ \mathfrak B$ has an element $y_D$ such that 
$ {\mathfrak B} \models \alpha  < y_D < \mu $ for every $ \alpha  < \mu $.


Since $ {\mathfrak B} $ is the Skolem hull
of $Id_D$ in ${\mathfrak C}$, we have $y_D =f(Id_D)$,
that is, $y_D =f_D$, for
some function $f: \lambda \to \lambda$ definable in ${\mathfrak A}$.
Since $f$ is definable, then also the following function $f'$
is definable:
\[
f'( \gamma )=
\begin{cases}
f( \gamma) & \textrm{ if } f(\gamma )< \mu \\
0 & \textrm{ if } f(\gamma ) \geq \mu \\
\end{cases}
\]
Since $ {\mathfrak B} \models  y_D < \mu $, then
$ \{ \gamma < \lambda | y( \gamma )< \mu \} \in D$.
Since
$y_D =f_D$,
$ \{ \gamma < \lambda | y( \gamma )=f( \gamma ) \} \in D$.
Hence,
$ \{ \gamma < \lambda | y( \gamma ) = f'( \gamma )\} \in D$,
being larger than the intersection of two sets in $D$.
Thus, 
$y_D =f'_D$.

Since $f': \lambda \to \mu$ and
$f'$ is definable  in ${\mathfrak A}$, then 
$f= f_ \beta $ for some $ \beta < \kappa $, thus
$y_D =(f_ \beta)_D $.


We need to show that $D'= f_ \beta (D)$ is uniform over $\mu$.
Indeed, for every $ \alpha_0 < \mu$,
and since $ {\mathfrak B} \models \alpha_0  < y_D $, then
$ \{ \gamma < \lambda | \alpha_0 < y( \gamma ) \} \in D $;
that is,
$ \{ \gamma < \lambda | \alpha_0 < f_ \beta ( \gamma ) \} \in D $, that is,
$ \{ \alpha < \mu| \alpha _0< \alpha \} \in D' $, and this implies
that $D'$ is uniform over $ \mu$, since $\mu$ is regular.   

(b) $\Rightarrow$ (a). Suppose we have functions
$ (f_ \beta ) _{ \beta < \kappa } $ as given by (b).

Expand 
$ \langle \lambda, <, \gamma \rangle _{ \gamma < \lambda }  $ to a model 
$ \mathfrak A$ 
by adding, for each $ \beta < \kappa $, a new function symbol representing
$f_ \beta $ (by abuse of notation, in what follows we shall
write $f_ \beta $ both for the function itself and
for the symbol that represents it).

Suppose that $\mathfrak B \equiv \mathfrak A$ and 
$ \mathfrak B$ has an element $x$ such that 
$ {\mathfrak B} \models \gamma < x $ for every $ \gamma < \lambda $.

For every formula $ \phi (z)$ with just one variable $z$ 
in the language of  
$ \mathfrak A$ let 
$E_ \phi = \{ \gamma < \lambda |  \mathfrak A \models \phi ( \gamma )\}  $.
Let $F= \{ E_ \phi | \mathfrak B \models \phi (x)\} $.
Since the intersection of any two members of $F$ is still in 
$F$, and $\emptyset \not\in F$, then $F$ can be extended to 
an ultrafilter $D$ on $ \lambda $. 

For every $ \gamma_0 < \lambda $,
consider the formula $ \phi (z) \equiv \gamma_0 < z$. We get
$E_ \phi = \{ \gamma < \lambda |  
\mathfrak A \models \gamma_0 < \gamma \}=
( \gamma_0 , \lambda )  $. On the other side, 
since $ {\mathfrak B} \models \gamma_0 < x $, then
by the definition of $F$ we have
$ E_ \phi = ( \gamma_0 , \lambda ) \in F \subseteq D$.
Thus, $D$ is uniform over $ \lambda $.  

By (b), 
$ f_ \beta (D) $ is uniform over  $ \mu$, for
some $ \beta < \kappa$. This means that
$( \alpha _0, \mu) \in f_ \beta (D)$, 
for every $ \alpha _0 < \mu$.
That is,
$ \{ \gamma < \lambda | \alpha _0 < f_ \beta ( \gamma )\} \in D $
for every $ \alpha _0 < \mu$.

For every $ \alpha _0 < \mu$,
consider the formula 
$ \psi (z) \equiv \alpha _0 < f_ \beta (z)$.
By the previous paragraph, 
$E_ \psi \in D$.
Notice that 
$E _{\neg \psi}  $ is the complement of 
$E_ \psi $ in $ \lambda $.
Since $D$ is proper,
and $ E_ \psi \in D$, then
$E _{\neg \psi} \not\in D $.
Since $D$ extends $F$,
and 
either $E_ \psi \in F$
or $E _{\neg \psi}  \in F$,
 we necessarily have
 $E_ \psi \in F$, that is, 
$\mathfrak B \models \psi (x)$, that is, 
$\mathfrak B \models \alpha _0 < f_ \beta (x)$.

Since $ \alpha _0 < \mu$ has been chosen
arbitrarily, we have that
$\mathfrak B \models \alpha _0 < f_ \beta (x)$
for every $ \alpha _0 < \mu$.
Moreover, since $ f_ \beta : \lambda \to \mu$,
and $ \mathfrak B \equiv \mathfrak A$, then 
$\mathfrak B \models f_ \beta (x) < \mu$.

Thus, we have proved that 
$ \mathfrak B$ has an element $y= f_ \beta (x)$ such that 
$ {\mathfrak B} \models \alpha  < y < \mu $ for every $ \alpha  < \mu $.

(b) $\Leftrightarrow$ (b$'$) follows from Lemma \ref{lem} below.

(b$'$) $ \Rightarrow $ (c).
Suppose that we have functions $ (f_ \beta ) _{ \beta < \kappa } $ as
given by (b$'$).
For $\alpha<\mu $
and $\beta<\kappa$, define 
$ B_{ \alpha , \beta } = f_ \beta ^{-1} ([0, \alpha ))  $. 

The family $ (B_{ \alpha , \beta }) _{ \alpha<\mu , \beta<\kappa}$
trivially satisfies Conditions (i) and (ii). Moreover, 
Condition (iii) is clearly equivalent to the condition
imposed on the $f_ \beta $'s in  (b$'$).

(c) $ \Rightarrow $ (b$'$).
Suppose we are given the family $ (B_{ \alpha , \beta }) _{ \alpha<\mu , \beta<\kappa}$
from (c). For $ \beta < \kappa $ and $ \gamma < \lambda $,  define
$f_ \beta ( \gamma )$ to be 
the smallest ordinal $ \alpha < \mu$  
such that $ \gamma \in B_{ \alpha , \beta }$ (such an $ \alpha $
exists because of (i)).

Because of Condition (ii),
we have that
$ B_{ \alpha , \beta } = f_ \beta ^{-1} ([0, \alpha ])  $,
for $\alpha<\mu $
and $\beta<\kappa$. Thus Condition (iii)
implies that 
for every function $g: \kappa \to \mu$ there exists some finite
set $F \subseteq \kappa $ such that 
$ \left| \bigcap _{\beta \in F} f_\beta  ^{-1}([0, g(\beta )]) \right| < \lambda $.
A fortiori, $ \left| \bigcap _{\beta \in F} f_\beta  ^{-1}([0, g(\beta ))) \right| < \lambda$, thus (b$'$) holds.

The equivalence of Conditions (c)-(e) has been proved in Part I \cite[Theorem 2]{parti}.
\end{proof}

\begin{lemma}\label{lem}
Suppose  that $ \lambda \geq \mu$  are infinite regular cardinals, 
and $\kappa  $ is a cardinal.  
Suppose that $ (f_ \beta ) _{ \beta < \kappa } $
is a given set of functions from $ \lambda $ to $\mu$.
Then the following are equivalent.

(a) Whenever $D$ is an ultrafilter
uniform over $ \lambda $ then there exists some $ \beta < \kappa $
such that $f_ \beta (D)$ is uniform over $ \mu$.

(b) For every function $g: \kappa \to \mu$ there exists some finite
set $F \subseteq \kappa $ such that 
$ \left| \bigcap _{\beta \in F} f_\beta  ^{-1}([0, g(\beta ))) \right| < \lambda $.
\end{lemma}

\begin{proof}
We show that the negation of (a) is equivalent to the negation of (b).

Indeed, (a) is false if and only if there exists an ultrafilter 
$D$ uniform over $\lambda $ such that for every $\beta<\kappa $ 
$f_\beta(D)$ is not uniform over $\mu$.
This means that for every $ \beta < \kappa $ there exists some
$g(\beta) < \mu$ such that 
$ [g(\beta) , \mu) \not\in f_ \beta (D)$,
that is, 
$ f_\beta ^{-1} ([g(\beta) , \mu)) \not\in D$,
that is, 
$ f_\beta ^{-1} ([0,g(\beta) )) \in D$.

Thus, there exists some $D$ which makes (a) false if and only if
there exists some function $g:\kappa\to \mu$
such that the set
$  \{ f_\beta ^{-1} ([0,g(\beta) )) | \beta < \kappa \} \cup \{ [ \gamma , \lambda )| \gamma<\lambda \}  $ has the finite intersection property.
Equivalently, there exists some function $g:\kappa\to \mu$
such that for every $F \subseteq \kappa $ 
the cardinality of $ \bigcap_{\beta\in F}  f_\beta ^{-1} ([0,g(\beta) )) $ 
is equal to $  \lambda  $ (since $ \lambda $ is regular).

This is exactly the negation of (b).
\end{proof}

\end{document}